\documentclass[11pt,a4paper]{article}
\usepackage{graphicx}
\usepackage{amsmath}

\newtheorem{theorem}{Theorem}

\newtheorem{corollary}[theorem]{Corollary}

\newtheorem{definition}[theorem]{Definition}

\newtheorem{lemma}[theorem]{Lemma}

\newtheorem{proposition}[theorem]{Proposition}

\newenvironment{proof}[1][Proof]{\textbf{#1.} }{\ \rule{0.5em}{0.5em}}

\begin{document}

\title{Stopping Games in Continuous Time\thanks{%
The results presented in this paper were proven while the authors
attended the workshop on ``Stochastic Methods in Decision and Game
Theory'', organized by Marco Scarsini in June 2002, Erice, Sicily,
Italy. The research of the second author was supported by the
Israel Science Foundation (grant No. 69/01-1).} }
\author{Rida Laraki \thanks{%
CNRS \textit{and} Laboratoire d'Econom\'{e}trie de l'Ecole Polytechnique, 1,
rue Descartes, 75005 Paris, France. Email: laraki@poly.polytechnique.fr} \
and Eilon Solan\thanks{%
MEDS Department, Kellogg School of Management, Northwestern
University, \emph{and} School of Mathematical Sciences, Tel Aviv
University, Tel Aviv 69978, Israel. Email: eilons@post.tau.ac.il,
e-solan@kellogg.northwestern.edu} } \maketitle

\begin{abstract}
We study two-player zero-sum stopping games in continuous time and infinite
horizon. We prove that the value in randomized stopping times exists as soon
as the payoff processes are right-continuous. In particular, as opposed to
existing literature, we do \emph{not} assume any conditions on the relations
between the payoff processes. We also show that both players have simple $%
\varepsilon $-optimal randomized stopping times; namely,
randomized stopping times which are small perturbations of
non-randomized stopping times.
\end{abstract}

\bigskip\bigskip

\noindent \textbf{Keywords: } Dynkin games, stopping games, optimal
stopping, stochastic analysis, continuous time.\bigskip

\newpage

\section{Introduction}

Stopping games in discrete time were introduced by Dynkin (1969) as a
variation of optimal stopping problems. In Dynkin's (1969) setup, two
players observe the realization of two discrete time processes $%
(x_{t},r_{t})_{t\in \mathbf{N}}$. Player 1 chooses a stopping time $\mu $
such that $\{\mu =t\}\subseteq \{r_{t}\geq 0\}$ for every $t\in {\mathbf{N}}$%
, and player 2 chooses a stopping time $\nu $ such that $\{\nu =t\}\subseteq
\{r_{t}<0\}$ for every $t\in {\mathbf{N}}$. Thus, players are not allowed to
stop simultaneously. Player 2 then pays player 1 the amount $x_{\min \{\mu
,\nu \}}\mathbf{1}_{\min \{\mu ,\nu \}<+\infty }$, where $\mathbf{1}$ is the
indicator function. This amount is a random variable. Denote the \emph{%
expected payoff} player 1 receives by
\begin{equation*}
\gamma (\mu ,\nu )={\mathbf{E}}[x_{\min \{\mu ,\nu \}}\mathbf{1}_{\min \{\mu
,\nu \}<+\infty }].
\end{equation*}
Dynkin (1969) proved that the game admits a value; that is,
\begin{equation*}
\sup_{\mu }\inf_{\nu }\gamma (\mu ,\nu )=\inf_{\nu }\sup_{\mu }\gamma (\mu
,\nu ).
\end{equation*}

Since then many authors generalized this basic result, both in discrete time
and in continuous time.

In \emph{discrete time}, Neveu (1975) allows the players to stop
simultaneously, that is, he introduces three uniformly integrable adapted
processes $(a_{t},b_{t},c_{t})_{t\in {\mathbf{N}}}$, the two players choose
stopping times $\mu $ and $\nu $ respectively, and the payoff player 2 pays
player 1 is
\begin{equation*}
a_{\mu }\mathbf{1}_{\{\mu <\nu \}}+b_{\nu }\mathbf{1}_{\{\mu >\nu \}}+c_{\mu
}\mathbf{1}_{\{\mu =\nu <+\infty \}}.
\end{equation*}
Neveu (1975) provides sufficient conditions for the existence of the value.
One of the conditions he imposes is the following:

\begin{itemize}
\item  \textbf{Condition C}: $c_{t}=a_{t}\leq b_{t}$ for every $t \geq 0$.
\end{itemize}

It is well known that in general the value need not exist when condition C
is not satisfied. Rosenberg et al.~(2001) allow the players to choose \emph{%
randomized} stopping times, and they prove, in discrete time again, the
existence of the value in randomized stopping times. This result was
recently generalized by Shmaya and Solan (2002) to the existence of an $%
\varepsilon $-equilibrium in the non-zero-sum problem.

Several authors, including Bismut (1979), Alario-Nazaret et al.~(1982) and
Lepeltier and Maingueneau (1984) studied the problem in \emph{continuous time}%
. That is, the processes $(a_{t},b_{t},c_{t})_{t\geq 0}$ are in continuous
time, and the stopping times the players choose are $[0,+\infty ]$-valued.
The literature provides sufficient conditions, that include condition C, for
the existence of the value in pure (i.e. non-randomized) stopping times.

Touzi and Vieille (2002) study the problem in continuous time, without
condition C, played on a bounded interval $[0,T]$; that is, players must
stop before or at time $T$. They prove that if $(a_{t})_{t\geq 0}$ and $%
(b_{t})_{t\geq 0}$ are semimartingales continuous at $T$, and if $c_{t}\leq
b_{t}$ for every $t\in \lbrack 0,T]$, then the game admits a value in
randomized stopping times.

In the present paper we prove that every stopping game in continuous time
where $\left( a_{t}\right) _{t\geq 0}$ and $\left( b_{t}\right) _{t\geq 0}$
are right-continuous, and $\left( c_{t}\right) _{t\geq 0}$ is progressively
measurable, admits a value in randomized stopping times. In addition, we
construct $\varepsilon $-optimal strategies which are as close as one wishes
to pure (non-randomized) stopping times; roughly speaking, there is a
stopping time $\mu $ such that for every $\delta $ sufficiently small there
is an $\varepsilon $-optimal strategy that stops with probability 1 between
times $\mu $ and $\mu +\delta $. Finally, we construct an $\varepsilon$%
-optimal strategy in the spirit of Dynkin (1969) and we extend the model by
introducing cumulative payoffs and final payoffs.

Stopping games in continuous time were applied in various contexts. The one
player stopping problem (the Snell envelope) is used in finance for the
pricing of the American option, see, e.g., Bensoussan (1984) and Karatzas
(1988). More recently Cvitanic and Karatzas (1996) used stopping games for
the study of backward stochastic differential equation with reflecting
barriers, and Ma and Cvitanic (2001) for the pricing of ``the American game
option''. Ghemawat and Nalebuff (1985) used stopping games to study
strategic exit from a shrinking market.

\section{Model, literature and main result}

A \emph{two-player zero-sum stopping game in continuous time} $\Gamma $ is
given by:

\begin{itemize}
\item  A probability space $(\Omega ,\mathcal{A},P)$: $\left( \Omega ,%
\mathcal{A}\right) $ is a measurable space and $P$ is a $\sigma$-additive
probability measure on $\left( \Omega ,\mathcal{A}\right) .$

\item  A filtration in continuous time $\mathcal{F}=(\mathcal{F}_{t})_{t\geq
0}$ satisfying ``the usual conditions''. That is, $\mathcal{F}$ is
right-continuous, and $\mathcal{F}_{0}$ contains all $P$-null sets: for
every $B\in \mathcal{A}$ with $P(B)=0$ and every $A\subset B,$ one has $A
\in \mathcal{F}_{0}$.

Denote $\mathcal{F}_{\infty }:=\vee _{t\geq 0}\mathcal{F}_{t}$. Assume
without loss of generality that $\mathcal{F}_{\infty }=\mathcal{A}.$ Hence $%
(\Omega ,\mathcal{A},P)$ is a complete probability space.

\item  Three uniformly bounded $\mathcal{F}$-adapted processes $%
(a_{t},b_{t},c_{t})_{t\geq 0}$.\footnote{%
As we argue below (see Section \ref{sec payoff}) our results hold for a
larger class of payoff processes, that contains the class of uniformly
integrable payoff processes.}
\end{itemize}

A \emph{pure strategy} of player 1 (resp.~player 2) is a $\mathcal{F}$%
-adapted stopping time $\mu $ (resp.~$\nu $). We allow players to never
stop, by choosing $\mu $ (or $\nu $) to be equal to $+\infty $.

The game proceeds as follows. Player 1 chooses a pure strategy $\mu $, and
player 2 chooses simultaneously and independently a pure strategy $\nu $.
Player 2 then pays player 1 the amount $a_{\mu }\mathbf{1}_{\{\mu <\nu
\}}+b_{\nu }\mathbf{1}_{\{\mu >\nu \}}+c_{\mu }\mathbf{1}_{\{\mu =\nu
<+\infty \}}$, which is a random variable. The \emph{expected payoff} that
correspond to a pair of pure strategies $(\mu ,\nu )$ is
\begin{equation*}
\gamma (\mu ,\nu )={\mathbf{E}}_{P}[a_{\mu }\mathbf{1}_{\{\mu <\nu
\}}+b_{\nu }\mathbf{1}_{\{\mu >\nu \}}+c_{\mu }\mathbf{1}_{\{\mu =\nu
<+\infty \}}].
\end{equation*}
Thus, if the game never stops, the payoff is $0$. This could be relaxed by
adding to the payoff a final payoff $\chi \mathbf{1}_{\mu=\nu=+\infty}$,
where $\chi$ is some $\mathcal{A}$-measurable function; see Section \ref
{generalization}. For a given stopping game $\Gamma $ we denote the expected
payoff by $\gamma _{\Gamma }(\mu ,\nu )$ when we want to emphasize the
dependency of the expected payoff on the game.

The quantity $\sup_{\mu }\inf_{\nu }\gamma (\mu ,\nu )$ is the maximal
amount that player 1 can guarantee to receive; that is, the best he can get
(in expectation) if player 2 knows the strategy chosen by player 1 before he
has to choose his own strategy. Similarly, by playing properly, player 2 can
guarantee to pay no more than $\inf_{\nu }\sup_{\mu }\gamma (\mu ,\nu )$.

\begin{definition}
If $\sup_{\mu }\inf_{\nu }\gamma (\mu ,\nu )=\inf_{\nu }\sup_{\mu }\gamma
(\mu ,\nu )$ then the common value is the \emph{value in pure strategies} of
the game, and is denoted by $v$. Any strategy $\mu $ for which $\inf_{\nu
}\gamma (\mu ,\nu )$ is within $\varepsilon $ of $v$ is \emph{an $%
\varepsilon $-optimal} strategy of player 1. $\varepsilon $-optimal
strategies of player 2 are defined analogously.
\end{definition}

Many authors provided sufficient conditions for the existence of the value
in pure strategies and $\varepsilon $-optimal pure strategies. The most
general set of sufficient conditions in continuous time was given by
Lepeltier and Maingueneau (1984, Corollary 12, Theorems 13 and 15).

\begin{theorem}[Lepeltier and Maingueneau, 1984]
\label{theorem lm} If (a) the processes $(a_{t},b_{t})_{t\geq 0}$ are
right-continuous, and (b) $a_{t}=c_{t}\leq b_{t}$ for every $t\geq 0,$ the
value exists and both players have pure $\varepsilon $-optimal strategies.
\end{theorem}

\textbf{Remark 1: } Lepeltier and Maingueneau (1984) require that the
processes $(a_{t},b_{t})_{t\geq 0}$ are optional; that is, measurable with
respect to the optional filtration. Recall that the optional filtration is
the one generated by all RCLL (right-continuous with left limit) processes.
Under the ``usual conditions'' it is also the filtration generated by all
right-continuous processes (see, e.g., Dellacherie and Meyer, 1975, \S IV-65).

\bigskip

Laraki (2000, Theorem 9.1) slightly extended this result by requiring that $%
c_{t}$ is in the convex hull of $a_{t}$ and $b_{t}$ ($c_{t}\in \mathrm{co}%
\{a_{t},b_{t}\}$) for every $t\geq 0$ instead of (b) of Theorem \ref{theorem
lm}.

The pure $\varepsilon $-optimal strategies that exist by Theorem \ref
{theorem lm} need not be finite. Indeed, the value of the game that is given
by $a_{t}=c_{t}=-1$ and $b_{t}=1$ for every $t\geq 0$ is 0, and the only $0$%
-optimal pure strategy of player 1 is $\mu =+\infty $. Moreover, if $\mu $
is an $\varepsilon $-optimal pure strategy of player 1 then $P(\mu <+\infty
)\leq \varepsilon $.

It is well known that in general the value in pure strategies need not
exist. Indeed, take $a_{t}=b_{t}=1$ and $c_{t}=0$ for every $t \geq 0$. Since $%
\gamma (\mu ,\mu )=0$ it follows that $\sup_{\mu }\inf_{\nu }\gamma (\mu
,\nu )=0$. For every stopping time $\nu $ define a stopping time $\mu _{\nu
} $ by
\begin{equation*}
\mu _{\nu }\left\{
\begin{array}{lll}
0 & \ \ \ \ \  & \nu >0 \\
1 &  & \nu =0
\end{array}
\right..
\end{equation*}
Since $\gamma (\mu _{\nu },\nu )=1$ for every $\nu $ it follows that $%
\inf_{\nu }\sup_{\mu }\gamma (\mu ,\nu )=1$, and the value in pure
strategies does not exist.

The difficulty with the last example is that
player 2, knowing the strategy of player 1, can stop exactly at the same
time as his opponent. The solution is to allow player 1 to choose his
stopping time randomly, thereby making the probability that the players stop
simultaneously vanish. Indeed, in the last example, if player 1 could have randomly chosen his
stopping time, say, uniformly in the interval $[0,1]$, then the game
terminates before time 1 with probability 1, and the probability of
simultaneous stopping is 0, whatever player 2 plays. In particular, such a
strategy guarantees player 1 payoff 1.

In the game theoretic literature, a standard and natural way to increase the
set of strategies is by allowing players to randomize. A \emph{mixed strategy%
} is a probability distribution over pure strategies. In general, this
allows to convexify the set of strategies, and makes the payoff function
bilinear. One can then apply a standard min-max theorem (e.g., Sion, 1958)
to prove the existence of the value in mixed strategies, provided some
regularity conditions hold (e.g., the space of mixed strategies is compact
and the payoff function continuous).

Three equivalent ways to randomize the set of pure strategies in our setup
are discussed in Touzi and Vieille (2002). We adopt the following definition
of mixed strategies due to Aumann (1964). It extends the probability space
to $([0,1]\times \lbrack 0,1]\times \Omega ,\mathcal{B\times B}\times
\mathcal{A},\lambda \otimes \lambda \otimes P)$, where $\mathcal{B}$ is the $%
\sigma $-algebra of Borel sets of $[0,1]$, and $\lambda $ is the Lebesgue
measure on $[0,1].$

\begin{definition}
A \emph{mixed strategy} of player 1 is a measurable function $\phi
:[0,1]\times \Omega \rightarrow \lbrack 0,+\infty ]$ such that for $\lambda $%
-almost every $r\in \lbrack 0,1]$, $\mu _{r}\left( \omega \right) :=\phi
(r,\omega )$ is a stopping time.
\end{definition}

The interpretation is the following: player 1 chooses randomly $r\in \lbrack
0,1]$, and then stops the game at time $\mu _{r}=\phi (r,\cdot )$. Mixed
strategies of player 2 are denoted by $\psi $, and, for every $s\in \lbrack
0,1]$, the $s$-section is denoted by $\nu _{s}:=\psi \left( s,\cdot \right) $%
.

The \emph{expected payoff} that corresponds to a pair of mixed strategies $%
(\phi ,\psi )$ is:
\begin{eqnarray}  \label{equ payoff}
\gamma (\phi ,\psi ) &=&\int_{[0,1]^{2}}\gamma (\mu _{r},\nu _{s})\ dr\ ds \\
&=&{\mathbf{E}}_{\lambda \otimes \lambda \otimes P}\left[a_{\mu _{r}}\mathbf{%
1}_{\{\mu _{r}<\nu _{s}\}}+b_{\nu _{s}}\mathbf{1}_{\{\mu _{r}>\nu
_{s}\}}+c_{\mu _{r}}\mathbf{1}_{\{\mu _{r}=\nu _{s}<+\infty \}}\right].
\notag
\end{eqnarray}

Though the payoff function given by (\ref{equ payoff}) is bilinear, without
strong assumptions on the data of the game the payoff function is not
continuous for the same topology which makes the strategy space compact.

\begin{definition}
If $\sup_{\phi }\inf_{\psi }\gamma (\phi ,\psi )=\inf_{\psi }\sup_{\phi
}\gamma (\phi ,\psi )$ then the common value is the \emph{value in mixed
strategies}, and it is denoted by $V$. Every strategy $\phi $ such that $%
\inf_{\psi }\gamma (\phi ,\psi )$ is within $\varepsilon $ of $V$ is \emph{$%
\varepsilon $-optimal} for player 1. $\varepsilon $-optimal strategies of
player 2 are defined analogously.
\end{definition}

Observe that $\sup_{\phi }\inf_{\psi }\gamma (\phi ,\psi )=\sup_{\phi
}\inf_{\nu }\gamma (\phi ,\nu )$, where $\nu$ ranges over all pure stopping times,
and that $\inf_{\psi }\sup_{\phi }\gamma
(\phi ,\psi )=\inf_{\psi }\sup_{\mu }\gamma (\mu ,\psi )$,
where $\mu$ ranges over all pure stopping times. Hence, to prove
the existence of the value, it suffices to show that $\sup_{\phi }\inf_{\nu
}\gamma (\phi ,\nu )=\inf_{\psi }\sup_{\mu }\gamma (\mu ,\psi )$. Moreover,
one always has $\sup_{\phi }\inf_{\psi }\gamma (\phi ,\psi )\leq \inf_{\psi
}\sup_{\phi }\gamma (\phi ,\psi )$.

Existence of the value in mixed strategies in stopping games with continuous
time was studied by Touzi and Vieille (2002), who proved the following.

Let $\Phi _{T}$ be the space of all mixed strategies $\phi $ such that $%
\lambda \otimes P(\mu _{r}\leq T)=1$, and let $\Psi _{T}$ be the space of
all mixed strategies $\psi $ such that $\lambda \otimes P(\nu _{s}\leq T)=1$.

\begin{theorem}[Touzi and Vieille, 2002]
\label{theorem tv} For every $T>0$, if (a) the processes $%
(a_{t},b_{t})_{t\geq 0}$ are semimartingales with trajectories continuous at
time $T$, (b) $c_{t}\leq b_{t}$ for every $t\geq 0$, and (c) the payoff
processes are uniformly integrable, then:
\begin{equation*}
\sup_{\phi \in \Phi _{T}}\inf_{\psi \in \Psi _{T}}\gamma (\phi ,\psi
)=\inf_{\psi \in \Psi _{T}}\sup_{\phi \in \Phi _{T}}\gamma (\phi ,\psi ).
\end{equation*}
\end{theorem}

Touzi and Vieille (2002) prove that under conditions (a) and (b) of Theorem
\ref{theorem tv} it is sufficient to restrict the players to certain
subclasses of mixed strategies. They then apply Sion's (1958) min-max
theorem to the restricted game.

\bigskip

\noindent \textbf{Remark 2:} By Dellacherie and Meyer (1980, \S VII-23),
under the ``usual conditions'', a semimartingle is always RCLL
(right-continuous with left limit).

\bigskip

One class of mixed strategies will play a special role along the paper.

\begin{definition}
Let $\delta > 0$. A mixed strategy $\phi$ is \emph{$\delta$-almost pure} if
there exists a stopping time $\mu$ and a set $A \in \mathcal{F}_\mu$ such
that for every $r \in [0,1]$, $\phi(r,\cdot) = \mu$ on $A$, and $%
\phi(r,\cdot) = \mu + r\delta$ on $A^c$.
\end{definition}

Recall that a process $\left( x_{t}\right) _{t\geq 0}$ is progressively
measurable if for every $t\geq 0$ the function $(s,\omega )\mapsto
x_{s}(\omega )$ from $\left[ 0,t\right] \times \Omega $ is measurable with
respect to $\mathcal{B}(\left[ 0,t\right] )\times \mathcal{F}_{t}$, where $%
\mathcal{B}(\left[ 0,t\right] )$ is the $\sigma $-algebra of Borel subsets
of $[0,t]$. Recall also that an optional process is progressively measurable
(see, e.g., Dellacherie and Meyer, 1975, \S IV-64).

The main result we present is the following.

\begin{theorem}
\label{theorem 1} If the processes $(a_{t})_{t\geq 0}$ and $(b_{t})_{t\geq
0} $ are right-continuous and if $(c_{t})_{t \geq 0}$ is progressively
measurable then the value in mixed strategies exists. Moreover, for every $%
\varepsilon >0$ there is $\delta _{0}\in (0,1)$ such that for every $\delta
\in (0,\delta _{0})$ both players have $\delta $-almost pure $\varepsilon $%
-optimal strategies.
\end{theorem}

Our proof heavily relies on the result of Lepeltier and Maingueneau (1984),
where they extend the discrete time variational approach of Neveu (1975) to
continuous time.

\subsection{On the payoff processes}

\label{sec payoff}

A $\mathcal{F}$-adapted process $x=\left( x_{t}\right) _{t\geq 0}$ is \emph{%
in the class $\mathcal{D}$} \ (see, e.g., Dellacherie and Meyer, 1980, \S
VI-20) if the set $\{x_{\sigma }\mathbf{1}_{\left\{ \sigma <+\infty \right\}
},\sigma \hbox{ is a }\mathcal{F}\hbox{-adapted stopping time}\}$ is
uniformly integrable (see, e.g., Dellacherie and Meyer, 1975, \S II-17). That
is, if for every bounded stopping time $\sigma $, $\mathbf{E}_{P}[|x_{\sigma
}|\mathbf{1}_{\left\{ x_{\sigma }\geq r\right\} }]$ converges uniformly to $%
0 $ as $r$ goes to $+\infty $.

This implies that the set $\{\mathbf{E}_{P}[|x_{\sigma }|\mathbf{1}_{\left\{
\sigma <+\infty \right\} }],\sigma \hbox{ is a }\mathcal{F}%
\hbox{-adapted
stopping time}\}$ is uniformly bounded (see, e.g., Dellacherie and Meyer,
1975, \S II-19). Observe that every uniformly bounded process, as well as
every uniformly integrable process, is in the class $\mathcal{D}$ (see, e.g.,
Dellacherie and Meyer, 1975, \S II-18)).

For a measurable process $\left( x_{t}\right) _{t\geq 0}$ and $r\geq 0$,
define the process $\left( x_{t}^{r}\right) _{t\geq 0}$ by:
\begin{equation*}
x_{t}^{r}(\omega ):=x_{t}(\omega )\mathbf{1}_{\{\left| x_{t}(\omega )\right|
\leq r\}}+r\mathbf{1}_{\{x_{t}(\omega )>r\}}-r\mathbf{1}_{\{x_{t}(\omega
)<-r\}}.
\end{equation*}
By Dellacherie and Meyer (1975, \S II-17), $x\in \mathcal{D}$ if and only if
for every $\varepsilon >0$ there exists $r>0$ such that for every stopping
time $\sigma $ one has $E\left[ \left| x_{\sigma }-x_{\sigma }^{r}\right|
\mathbf{1}_{\left\{ \sigma <+\infty \right\} }\right] <\varepsilon $. The
process $\left( x_{t}^{r}\right) _{t\geq 0}$ is uniformly bounded by $r$,
and, in addition, if $\left( x_{t}\right) _{t\geq 0}$ is right-continuous or
$\mathcal{F}$-adapted, so is $\left( x_{t}^{r}\right) _{t\geq 0}.$

If the payoff processes $(a_{t})_{t\geq 0}$, $(b_{t})_{t\geq 0}$ and $%
(c_{t})_{t\geq 0}$ are not necessarily bounded but are in the class $%
\mathcal{D}$, then for every $\varepsilon >0$ there exists $r>0$ such that
\begin{equation*}
{\mathbf{E}}_{P}\left[ \left( \left| a_{\sigma }-a_{\sigma }^{r}\right|
+\left| b_{\sigma }-b_{\sigma }^{r}\right| +\left| c_{\sigma }-c_{\sigma
}^{r}\right| \right) \mathbf{1}_{\left\{ \sigma <+\infty \right\} }\right]
<\varepsilon .
\end{equation*}
Hence, if the game $\Gamma =(\Omega ,\mathcal{A},P;\mathcal{F}%
,(a_{t},b_{t},c_{t})_{t\geq 0})$ satisfies the assumptions of Theorem \ref{theorem 1}, it
admits a value. Moreover, every $\varepsilon $-optimal strategy in $\Gamma
^{r}:=(\Omega ,\mathcal{A},P;\mathcal{F},(a_{t}^{r},b_{t}^{r},c_{t}^{r})_{t%
\geq 0})$ is $2\varepsilon $-optimal in $\Gamma $.

In particular, all the existence results that are proved for uniformly
bounded payoff processes (Lepeltier and Maingueneau (1984)) or uniformly
integrable payoff processes (Touzi and Vieille (2002)) are valid for payoff
processes in the class $\mathcal{D}$ as well.

\section{Proof}

In the present section the main result of the paper is proven. From now on
we fix a stopping game $\Gamma $ such that $(a_{t},b_{t})_{t\geq 0}$ are
right-continuous and $\left( c_{t}\right) _{t\geq 0}$ is progressively
measurable.

\subsection{Preliminaries}

The following Lemma will be used in the sequel.

\begin{lemma}
\label{lemma 1} For every $\mathcal{F}$-adapted stopping time $\tau $ and
every $\varepsilon >0$ there is $\delta >0$ such that $p\left(\{|a_{t}-a_{%
\tau }|<\varepsilon \ \ \forall t\in [ \tau ,\tau +\delta
]\}\right)>1-\varepsilon $.
\end{lemma}

A similar statement holds when one replaces the process $(a_{t})_{t\geq 0}$
by the process $(b_{t})_{t\geq 0}$.

\begin{proof}
Since $(a_{t})_{t\geq 0}$ is right-continuous, it is progressively
measurable (see, e.g., Dellacherie and Meyer, 1975, \S IV-15)$.$

Let $\delta _{\tau (\omega )}(\omega )=\inf \{s\geq \tau (\omega ):$ $%
|a_{s}(\omega )-a_{\tau }(\omega )|\geq \varepsilon \}.$ The progressive
measurability of $(a_{t})_{t\geq 0}$ implies that $\delta _{\tau (\omega
)}(\omega )$ is measurable with respect to $\mathcal{F}_{\infty }$ (see,
e.g., Dellacherie and Meyer, 1975, \S III-44)$.$

The right-continuity of $(a_{t})_{t\geq 0}$ implies that $P(\{\omega :\delta
_{\tau (\omega )}(\omega )>0\})=1$. Since $P$ is $\sigma $-additive and $%
\{\omega :\delta _{\tau (\omega )}(\omega )>0\}=\cup _{n>0}\{\omega :\delta
_{\tau (\omega )}(\omega )>\frac{1}{n}\}$, the lemma follows by choosing $%
\delta >0$ sufficiently small so that $P(\{\omega :\delta _{\tau (\omega
)}(\omega )>\delta \})>1-\varepsilon $.
\end{proof}

\bigskip

By Lemma \ref{lemma 1}, and since the payoff processes are uniformly bounded,
one obtains the following.

\begin{corollary}
\label{corollary 1} Let a stopping time $\tau $ and $\varepsilon >0$ be
given. There exists $\delta >0$ sufficiently small such that for every $\mathcal{F}_{\tau }$-measurable set $%
A \subseteq \{\tau < +\infty\}$, and every stopping time $\mu $ that satisfies $%
\tau \leq \mu \leq \tau +\delta $,
\begin{equation*}
\left| {\mathbf{E}}_{P}[a_{\mu }\mathbf{1}_{A}]-{\mathbf{E}}_{P}[a_{\tau }%
\mathbf{1}_{A}]\right| \leq 2\varepsilon .
\end{equation*}
\end{corollary}

\subsection{The case $a_t \leq b_t$ for every $t \geq 0$}

\begin{definition}
Let $\delta > 0$. A mixed strategy $\phi$ is \emph{$\delta$-pure} if there
exists a stopping time $\mu$ such that
\begin{equation}  \label{equ delta pure}
\phi(r,\cdot) = \mu + r\delta \ \ \ \ \ \forall r \in [0,1].
\end{equation}
\end{definition}

Observe that a $\delta$-pure mixed strategy is in particular $\delta$-almost
pure. When $\mu$ is a stopping time, we sometime denote the $\delta$-pure
mixed strategy defined in (\ref{equ delta pure}) simply by $\mu + r\delta$.

In this section we prove the following result: when $a_{t}\leq b_{t}$ for
every $t\geq 0$ the value in mixed strategies exists, it is independent of $%
(c_{t})_{t\geq 0}$, and both players have $\delta$-pure $\epsilon$-optimal
strategies, provided $\delta$ is sufficiently small.

The idea is the following. Assume player 1 decides to stop at time $t$. If $%
c_{t}\leq a_{t}$, player 1 wants to mask the exact time in which he stops,
so that player 2 cannot stop at the same time. Since payoffs are
right-continuous, he can stop randomly in a small interval after time $t$.
If $c_{t}>a_{t}$, player 2 prefers that player 1 stops alone at
time $t$ rather than to stop simultaneously with player 1 at time $t$.

\begin{proposition}
\label{proposition 1} If $a_{t}\leq b_{t}$ for every $t\geq 0$ then the
value in mixed strategies exists. Moreover, the value is independent of the
process $(c_{t})_{t\geq 0}$, and for every $\varepsilon >0$ there is $\delta
_{0}>0$ such that for every $\delta \in (0,\delta _{0})$ both players have $%
\delta $-pure $\varepsilon $-optimal strategies. If $a_{t}\leq c_{t}\leq
b_{t}$ for every $t\geq 0$ then the value in pure strategies exists, and
there are $\epsilon$-optimal strategies that are independent of $%
(c_{t})_{t\geq 0}.$
\end{proposition}

\begin{proof}
Consider an auxiliary stopping game $\Gamma ^{\ast }=(\Omega ,\mathcal{A},P;%
\mathcal{F},(a_{t}^{\ast },b_{t}^{\ast },c_{t}^{\ast })_{t \geq 0})$, where $%
a_{t}^{\ast }=a_{t}$ and $b_{t}^{\ast }=c_{t}^{\ast}=b_{t}$ for every $t
\geq 0$.

By Theorem \ref{theorem lm} the game $\Gamma ^{\ast }$ has a value in pure
strategies $v^{\ast }$. We will prove that $v^{\ast }$ is the value in mixed
strategies of the original game. Since $\Gamma^{\ast }$ does not depend on the
process $(c_{t})_{t\geq 0}$, the second claim in the proposition will follow.

Fix $\varepsilon > 0$. Let $\mu$ be an $\varepsilon$-optimal strategy of
player 1 in $\Gamma^*$. In particular, $\inf_{\nu }\gamma _{\Gamma ^{\ast
}}(\mu,\nu )\geq v^*-\varepsilon $.

We now construct a mixed strategy $\phi $ that satisfies $\inf_{\nu }\gamma
_{\Gamma }(\phi,\nu )\geq v^{\ast }-5\varepsilon $. By Lemma \ref{lemma 1}
there is $\delta >0$ such that $p(\{|a_{t}-a_{\mu }|<\varepsilon \ \ \ \
\forall t\in \lbrack \mu ,\mu +\delta ]\})>1-\varepsilon $. Define a $\delta
$-pure mixed strategy $\phi $ by
\begin{equation*}
\phi (r,\cdot )=\mu +r\delta \ \ \ \ \ \forall r\in \lbrack 0,1].
\end{equation*}

Let $\nu$ be any stopping time. Since $\mu$ is $\varepsilon$-optimal in $%
\Gamma^*$, by the definition of $\Gamma^*$, and since $\lambda \otimes P(\mu
+ r \delta = \nu) = 0$,
\begin{eqnarray}  \label{aux1}
v^* - \varepsilon &\leq& \gamma_{\Gamma^*}(\mu,\nu)  \notag \\
&=& {\mathbf{E}}_P[ a_\mu {\mathbf{1}}_{\{\mu < \nu\}} + b_{\nu} {\mathbf{1}}%
_{\{\mu \geq \nu\}}] \\
&=& {\mathbf{E}}_{\lambda \otimes P}[ a_\mu {\mathbf{1}}_{\{\mu+r\delta <
\nu\}} + a_\mu {\mathbf{1}}_{\{\mu < \nu < \mu+r\delta\}} + b_{\nu} {\mathbf{%
1}}_{\{\mu \geq \nu\}}].  \notag
\end{eqnarray}

Since $\lambda \otimes P(\mu +r\delta =\nu )=0$ and $\left( c_{t}\right)
_{t\geq 0}$ is progressively measurable,
\begin{eqnarray}
\gamma _{\Gamma }(\phi ,\nu ) &=&{\mathbf{E}}_{\lambda \otimes P}\left[
a_{\mu +r\delta }\mathbf{1}_{\{\mu +r\delta <\nu \}}+b_{\nu }\mathbf{1}%
_{\{\mu +r\delta >\nu \}}+c_{\nu }\mathbf{1}_{\{\mu +r\delta =\nu <+\infty
\}}\right]  \notag  \label{aux2} \\
&=&{\mathbf{E}}_{\lambda \otimes P}\left[ a_{\mu +r\delta }\mathbf{1}_{\{\mu
+r\delta <\nu \}}+b_{\nu }\mathbf{1}_{\{\mu +r\delta >\nu \}}\right] \\
&=&{\mathbf{E}}_{\lambda \otimes P}\left[ a_{\mu +r\delta }\mathbf{1}_{\{\mu
+r\delta <\nu \}}+b_{\nu }\mathbf{1}_{\{\mu <\nu <\mu +r\delta \}}+b_{\nu }%
\mathbf{1}_{\{\mu \geq \nu \}}\right] .  \notag
\end{eqnarray}
By Corollary \ref{corollary 1}, and since $a_{t}\leq b_{t}$ for every $t\geq
0$,
\begin{equation}
{\mathbf{E}}_{\lambda \otimes P}[a_{\mu }\mathbf{1}_{\{\mu <\nu <\mu
+r\delta \}}]\leq {\mathbf{E}}_{\lambda \otimes P}[a_{\nu }\mathbf{1}_{\{\mu
<\nu <\mu +r\delta \}}]+2\varepsilon \leq {\mathbf{E}}_{\lambda \otimes
P}[b_{\nu }\mathbf{1}_{\{\mu <\nu <\mu +r\delta \}}]+2\varepsilon .
\label{aux3}
\end{equation}
Corollary \ref{corollary 1} implies in addition that
\begin{equation}
{\mathbf{E}}_{\lambda \otimes P}[a_{\mu }\mathbf{1}_{\{\mu +r\delta <\nu
\}}]\leq {\mathbf{E}}_{\lambda \otimes P}[a_{\mu +r\delta }\mathbf{1}_{\{\mu
+r\delta <\nu \}}]+2\varepsilon .  \label{aux4}
\end{equation}
By (\ref{aux1})-(\ref{aux4}),
\begin{equation*}
v^{\ast }-\varepsilon \leq \gamma _{\Gamma ^{\ast }}(\mu ,\nu )\leq \gamma
_{\Gamma }(\phi ,\nu )+4\varepsilon .
\end{equation*}
Since $\nu $ is arbitrary, $\inf_{\nu }\gamma _{\Gamma }(\phi ,\nu )\geq
v^{\ast }-5\varepsilon $.

Consider an auxiliary stopping game $\Gamma ^{\ast \ast }=(\Omega ,\mathcal{A%
},P;\mathcal{F},(a_{t}^{\ast \ast },b_{t}^{\ast \ast },c_{t}^{\ast \ast
})_{t\geq 0})$, where $a_{t}^{\ast \ast }=c_{t}^{\ast \ast }=a_{t}$ and $%
b_{t}^{\ast \ast }=b_{t}$ for every $t\geq 0$.

A symmetric argument to the one provided above proves that the game $\Gamma
^{\ast \ast }$ has a value $v^{\ast \ast }$, and that player 2 has a mixed
strategy $\psi $ which satisfies $\sup_{\mu }\gamma _{\Gamma }(\mu ,\psi
)\leq v^{\ast \ast }+5\varepsilon $.

Since $c_{t}^{\ast \ast }=a_{t}\leq b_{t}=c_{t}^{\ast }$ for every $t\geq 0$%
, $v^{\ast \ast }\leq v^{\ast }$. Since $\sup_{\mu }\gamma _{\Gamma }(\mu
,\psi )\geq \gamma _{\Gamma }(\phi ,\psi )\geq \inf_{\nu }\gamma _{\Gamma
}(\phi ,\nu )$,
\begin{equation*}
v^* \geq v^{\ast \ast }\geq \sup_{\mu }\gamma _{\Gamma }(\mu ,\psi
)-5\varepsilon \geq \inf_{\nu }\gamma _{\Gamma }(\phi ,\nu )-5\varepsilon
\geq v^{\ast }-10\varepsilon .
\end{equation*}
Since $\varepsilon $ is arbitrary, $v^{\ast }=v^{\ast \ast }$, so that $%
v^{\ast }$ is the value in mixed strategies of $\Gamma $, and $\phi $ and $%
\psi $ are $5\varepsilon $-optimal mixed strategies of the two players.

If $a_{t}\leq c_{t}\leq b_{t}$ for every $t\geq 0$ then $\gamma _{\Gamma
^{\ast \ast }}(\mu ,\nu )\leq \gamma _{\Gamma }(\mu ,\nu )\leq \gamma
_{\Gamma ^{\ast }}(\mu ,\nu )$ for every pair of pure strategies $(\mu ,\nu
) $. Hence
\begin{eqnarray*}
v^{\ast \ast } &=&\sup_{\mu }\inf_{\nu }\gamma _{\Gamma ^{\ast \ast }}(\mu
,\nu )\leq \sup_{\mu }\inf_{\nu }\gamma _{\Gamma }(\mu ,\nu ) \\
&\leq &\inf_{\nu }\sup_{\mu }\gamma _{\Gamma }(\mu ,\nu )\leq \inf_{\nu
}\sup_{\mu }\gamma _{\Gamma ^{\ast }}(\mu ,\nu )=v^{\ast }=v^{\ast \ast }.
\end{eqnarray*}
Thus $\sup_{\mu }\inf_{\nu }\gamma _{\Gamma }(\mu ,\nu )=\inf_{\nu
}\sup_{\mu }\gamma _{\Gamma }(\mu ,\nu )$ : the value in pure strategies
exists. Moreover, any $\epsilon $-optimal strategy of player 1 (resp.~player
2) in $\Gamma ^{\ast }$ (resp.~$\Gamma ^{\ast \ast }$) is also $\epsilon $%
-optimal in $\Gamma $. In particular, if $a_{t}\leq c_{t}\leq b_{t}$ for
every $t\geq 0,$ both players have $\epsilon $-optimal strategies that are
independent of $(c_{t})_{t\geq 0}$.
\end{proof}

\subsection{Proof of Theorem \ref{theorem 1}}

Define a stopping time $\tau $ by
\begin{equation*}
\tau =\inf \{t\geq 0,a_{t}\geq b_{t}\},
\end{equation*}
where the infimum of an empty set is $+\infty $. Since $(a_{t}-b_{t})_{t\geq
0}$ is progressively measurable with respect to $\left( \mathcal{F}%
_{t}\right) _{t\geq 0}$, $\tau $ is an $\mathcal{F}$-adapted stopping time
(see, e.g., Dellacherie and Meyer, 1975, \S IV-50).

The idea is the following. We show that it is optimal for both players to
stop at or around time $\tau $ (provided the game does not stop before time $%
\tau ).$ Hence the problem reduces to the game between times $0$ and $\tau $%
. Since for $t\in \lbrack 0,\tau \lbrack $, $a_{t}\leq b_{t},$ Proposition
\ref{proposition 1} can be applied.

\bigskip 

The following notation will be useful in the sequel. For a pair of pure
strategies $(\mu ,\nu )$, and a set $A\in \mathcal{A}$, we define
\begin{equation*}
\gamma _{\Gamma }(\mu ,\nu ;A)={\mathbf{E}}_{P}[\mathbf{1}_{A}(a_{\mu }%
\mathbf{1}_{\{\mu <\nu \}}+b_{\mu }\mathbf{1}_{\{\mu >\nu \}}+c_{\mu }%
\mathbf{1}_{\{\mu =\nu <+\infty \}})].
\end{equation*}
This is the expected payoff restricted to $A$. For a pair of mixed
strategies $(\phi ,\psi )$ we define
\begin{equation*}
\gamma _{\Gamma }(\phi ,\psi ;A)=\int_{[0,1]^{2}}\gamma _{\Gamma }(\mu
_{r},\nu _{s};A)dr\ ds,
\end{equation*}
where $\mu _{r}$ and $\nu _{s}$ are the sections of $\phi $ and $\psi $
respectively.

Set
\begin{eqnarray*}
A_{0} &=&\{\tau =+\infty \}, \\
A_{1} &=&\{\tau <+\infty \}\cap \{c_{\tau }\geq a_{\tau }\geq b_{\tau }\}, \\
A_{2} &=&\{\tau <+\infty \}\cap \{a_{\tau }>c_{\tau }\geq b_{\tau }\},%
\hbox{
and} \\
A_{3} &=&\{\tau <+\infty \}\cap \{a_{\tau }\geq b_{\tau }>c_{\tau }\}.
\end{eqnarray*}
Observe that $(A_{0},A_{1},A_{2},A_{3})$ is an $\mathcal{F}_{\tau }$%
-measurable partition of $\Omega $.

Define a $\mathcal{F}_{\tau }$-measurable function $w$ by
\begin{equation*}
w=a_{\tau }\mathbf{1}_{A_{1}}+c_{\tau }\mathbf{1}_{A_{2}}+b_{\tau }\mathbf{1}%
_{A_{3}}.
\end{equation*}

Define a stopping game $\Gamma ^{\ast }=(\Omega ,\mathcal{A},P,(\mathcal{F}%
_{t})_{t\geq 0},(a_{t}^{\ast },b_{t}^{\ast },c_{t}^{\ast })_{t\geq 0})$ by:
\begin{equation*}
a_{t}^{\ast }=\left\{
\begin{array}{lll}
a_{t} & \  & t<\tau \\
w &  & t\geq \tau
\end{array}
\right. ,\ \ \ \ \ b_{t}^{\ast }=\left\{
\begin{array}{lll}
b_{t} & \  & t<\tau \\
w &  & t\geq \tau
\end{array}
\right. ,\ \ \ \ \ c_{t}^{\ast }=\left\{
\begin{array}{lll}
c_{t} & \  & t<\tau \\
w &  & t\geq \tau
\end{array}
\right. .
\end{equation*}
That is, the payoff is set to $w$ at and after time $\tau $.

The game $\Gamma ^{* }$ satisfies the assumptions of Proposition \ref
{proposition 1}, hence it has a value in mixed strategies $V$. Moreover, for
every $\varepsilon > 0$ both players have $\delta$-pure $\varepsilon $%
-optimal strategies, provided $\delta > 0$ is sufficiently small.

We now prove that $V$ is the value of the game $\Gamma$ as well. Fix $%
\varepsilon >0$. We only show that player 1 has a mixed strategy $\phi $
such that $\inf_{\nu }\gamma _{\Gamma }(\phi ,\nu )\geq V-7\varepsilon $. An
analogous argument shows that player 2 has a mixed strategy $\psi $ such
that $\sup_{\mu }\gamma _{\Gamma }(\mu ,\psi )\leq V+7\varepsilon $. Since $%
\varepsilon $ is arbitrary, $V$ is indeed the value in mixed strategies of $%
\Gamma $.

Assume $\delta$ is sufficiently small so that the following conditions hold
(by the proofs of Lemma \ref{lemma 1} and Proposition \ref{proposition 1}
such $\delta$ exists).

\begin{enumerate}
\item[(\textbf{C1})]  Player 1 has a $\delta$-pure $\varepsilon$-optimal strategy $%
\phi^* = \mu + r\delta$ in $\Gamma^*$.

\item[(\textbf{C2})]  $P(\{\mu +\delta <\tau \})\geq P(\{\mu <\tau \})-\varepsilon /M$%
, where $M \in ]0,+\infty[$ is a uniform bound of the payoff processes.

\item[(\textbf{C3})]  $P(\{|a_{t}-a_{\tau }|<\varepsilon ,|b_{t}-b_{\tau
}|<\varepsilon \ \ \ \forall t\in \lbrack \tau ,\tau +\delta
]\})>1-\varepsilon $.
\end{enumerate}

We now claim that one can choose $\mu $ so that $\mu \leq \tau $. Indeed,
assume that $P(\{\mu >\tau \})>0$. The set $\{\mu >\tau \}$ is $\mathcal{F}%
_{\tau }$-measurable. Define a stopping time $\mu ^{\prime }=\min \{\mu
,\tau \}$. We will prove that the $\delta $-pure strategy $\phi ^{\prime
}=\mu ^{\prime }+r\delta $ is also $\varepsilon $-optimal in $\Gamma ^{\ast
} $, which establishes the claim. Given a stopping time $\nu $ define a
stopping time $\nu ^{\prime }$ as follows: $\nu ^{\prime }=\tau $ over $%
\{\mu >\tau \}$, and $\nu ^{\prime }=\nu $ otherwise. Then
\begin{equation*}
V-\varepsilon \leq \gamma _{\Gamma ^{\ast }}(\mu +r\delta ,\nu ^{\prime
})=\gamma _{\Gamma ^{\ast }}(\mu +r\delta ,\nu ^{\prime };\{\mu >\tau
\})+\gamma _{\Gamma ^{\ast }}(\mu +r\delta ,\nu ^{\prime };\{\mu \leq \tau
\}).
\end{equation*}
However, $\gamma _{\Gamma ^{\ast }}(\mu +r\delta ,\nu ^{\prime };\{\mu >\tau
\})={\mathbf{E}}_{\lambda \otimes P}[w\mathbf{1}_{\{\mu > \tau\}}]=\gamma
_{\Gamma ^{\ast }}(\mu ^{\prime }+r\delta ,\nu ;\{\mu >\tau \})$, and since $%
\mu =\mu ^{\prime }$ and $\nu =\nu ^{\prime }$ over $\{\mu \leq \tau \}$, $%
\gamma _{\Gamma ^{\ast }}(\mu +r\delta ,\nu ^{\prime };\{\mu \leq \tau
\})=\gamma _{\Gamma ^{\ast }}(\mu ^{\prime }+r\delta ,\nu ;\{\mu \leq \tau
\} $. Therefore
\begin{equation*}
\gamma _{\Gamma ^{\ast }}(\mu ^{\prime }+r\delta ,\nu )=\gamma _{\Gamma
^{\ast }}(\mu +r\delta ,\nu ^{\prime })\geq V-\varepsilon .
\end{equation*}
Since $\nu $ is arbitrary, $\mu ^{\prime }+r\delta $ is $\varepsilon $%
-optimal, as desired.

Define a mixed strategy $\phi $ as follows.
\begin{equation*}
\phi (r,\cdot )=\left\{
\begin{array}{lll}
\mu +r\delta & \ \ \ \ \  & \{\mu <\tau \}\cup A_{0}, \\
\tau &  & \{\mu =\tau\}\cap \left( A_{1}\cup A_{2}\right), \\
\mu +r\delta &  & \{\mu =\tau \}\cap A_{3}.
\end{array}
\right.
\end{equation*}
Observe that $\phi$ is $\delta$-almost pure.

The mixed strategies $\phi$ and $\phi^*$ differ only over the set $\{\mu
=\tau\}\cap \left( A_{1}\cup A_{2}\right)$. Since over this set the payoff
in $\Gamma^*$ is $w$ provided the game terminates after time $\tau$,
whatever the players play, $\phi$ is an $\varepsilon$-optimal mixed strategy
in $\Gamma ^{* }$.

Let $\nu $ be an arbitrary pure strategy of player 2. Define a partition $%
(B_{0},B_{1},B_{2})$ of $[0,1]\times \Omega $ by
\begin{eqnarray*}
B_{0} &=&\{\mu +\delta <\tau \}\cup \{\nu <\tau \}, \\
B_{1} &=&\{\mu <\tau <\mu +\delta \}\cap \{\nu \geq \tau \},\hbox{ and} \\
B_{2} &=&\left( \{\mu =\tau \hbox{ or }\mu =+\infty \}\right) \cap \{\nu
\geq \tau \}.
\end{eqnarray*}

Over $B_{0}$ the game terminates before time $\tau $ under $(\phi ,\nu )$.
In particular,
\begin{equation}
\gamma _{\Gamma }(\phi ,\nu ;B_{0})=\gamma _{\Gamma ^{\ast }}(\phi ,\nu
;B_{0}).  \label{final1}
\end{equation}
By (\textbf{C2}) $\lambda \otimes P(B_{1})<\varepsilon/M$, so that
\begin{equation}
\gamma _{\Gamma }(\phi ,\nu ;B_{1})\geq \gamma _{\Gamma ^{\ast }}(\phi ,\nu
;B_{1})-2\varepsilon .  \label{final2}
\end{equation}
Over $B_{2}\cap A_{0}$ the game never terminates under $(\phi ,\nu )$, so
that
\begin{equation}
\gamma _{\Gamma }(\phi ,\nu ;B_{2}\cap A_{0})=\gamma _{\Gamma ^{\ast }}(\phi
,\nu ;B_{2}\cap A_{0})=0.  \label{final3}
\end{equation}
Over $A_{1}\cup A_{2}$, $\min \{a_{\tau },c_{\tau }\}\geq w$, so that
\begin{eqnarray}
\gamma _{\Gamma }(\phi ,\nu ;B_{2}\cap (A_{1}\cup A_{2})) &=&{\mathbf{E}}%
_{\lambda \otimes P}[\mathbf{1}_{B_{2}\cap (A_{1}\cup A_{2})}(a_{\tau }%
\mathbf{1}_{\{\tau <\nu \}}+c_{\tau }\mathbf{1}_{\{\tau =\nu \}})]  \notag
\label{final4} \\
&\geq &{\mathbf{E}}_{\lambda \otimes P}[w\mathbf{1}_{\{\tau \leq \nu \}\cap
B_{2}\cap (A_{1}\cup A_{2})}] \\
&=&\gamma _{\Gamma ^{\ast }}(\phi ,\nu ;B_{2}\cap (A_{1}\cup A_{2})).  \notag
\end{eqnarray}
Finally, since $\lambda \otimes P(\{\mu +r\delta =\nu \})=0$, by Corollary
\ref{corollary 1}, since $\left( c_{t}\right) _{t\geq 0}$ is progressively
measurable, and since $a_{\tau }\geq b_{\tau }=w$ over $A_{3}$,
\begin{eqnarray}
\gamma _{\Gamma }(\phi ,\nu ;B_{2}\cap A_{3}) &=&{\mathbf{E}}_{\lambda
\otimes P}[\mathbf{1}_{B_{2}\cap A_{3}}(a_{\mu +r\delta }\mathbf{1}_{\{\mu
+r\delta <\nu \}}+b_{\nu }\mathbf{1}_{\{\mu +r\delta >\nu \}}+c_{\nu }%
\mathbf{1}_{\{\mu +r\delta =\nu \}})]  \notag  \label{final5} \\
&=&{\mathbf{E}}_{\lambda \otimes P}[\mathbf{1}_{B_{2}\cap A_{3}}(a_{\mu
+r\delta }\mathbf{1}_{\{\mu +r\delta <\nu \}}+b_{\nu }\mathbf{1}_{\{\mu
+r\delta >\nu \}})]  \notag \\
&\geq &{\mathbf{E}}_{\lambda \otimes P}[\mathbf{1}_{B_{2}\cap A_{3}}(a_{\tau
}\mathbf{1}_{\{\mu +r\delta <\nu \}}+b_{\tau }\mathbf{1}_{\{\mu +r\delta
>\nu \}})]-4\varepsilon \\
&\geq &{\mathbf{E}}_{\lambda \otimes P}[w\mathbf{1}_{B_{2}\cap
A_{3}}]-4\varepsilon  \notag \\
&=&\gamma _{\Gamma ^{\ast }}(\phi ,\nu ;B_{2}\cap A_{3}).  \notag
\end{eqnarray}

Summing Eqs. (\ref{final1})-(\ref{final5}), and using the $\varepsilon $%
-optimality of $\phi ^{\ast }$ in $\Gamma ^{\ast }$, gives us
\begin{equation*}
V-\varepsilon \leq \gamma _{\Gamma ^{\ast }}(\phi ,\nu )\leq \gamma _{\Gamma
}(\phi ,\nu )+6\varepsilon ,
\end{equation*}
as desired.

\section{Extensions}

\label{generalization}

In the present section we construct specific $\varepsilon $-optimal
strategies in the spirit of Dynkin (1969) or Rosenberg et al.~(2001), and we
give conditions for the existence of the value in pure strategies. We then
provide two extensions to the basic model.

\subsection{Construction of an $\protect\varepsilon$-optimal strategy}

Let $\Gamma =(\Omega ,\mathcal{A},P,(\mathcal{F}_{t})_{t\geq
0},(a_{t},b_{t},c_{t})_{t\geq 0})$ satisfy the conditions of Theorem \ref{theorem 1}.

Define $\tau ,$ $(A_{0},A_{1},A_{2},A_{3})$, $w$ and $\Gamma ^{\ast
}=(\Omega ,\mathcal{A},P,(\mathcal{F}_{t})_{t\geq 0},(a_{t}^{\ast
},b_{t}^{\ast },c_{t}^{\ast })_{t\geq 0})$ as in the proof of Theorem \ref{theorem 1}.

For any stopping time $\sigma $ let $\Gamma _{\sigma }^{\ast }=(\Omega ,%
\mathcal{A},P,(\mathcal{F}_{t})_{t\geq 0},(a_{t}^{\ast },b_{t}^{\ast
},c_{t}^{\ast })_{t\geq 0})$ be the game starting at time $\sigma $; that
is, players are restricted to choose strategies that stop with probability 1
at or after time $\sigma $.

Lepeltier and Mainguenau (1984, Theorem 13) and Proposition \ref{proposition
1} show that this game has a value in mixed strategies $X_{\sigma}^{\ast }$.
Moreover, the value is independent of $\left( c_{t}^{\ast }\right) _{t\geq
0} $.

Using a general result of Dellacherie and Lenglart (1982), Lepeltier and
Mainguenau (1984, Theorem 7) show the existence of a right-continuous
process $\left( V_{t}^{\ast }\right) _{t\geq 0}$ such that $V_{\sigma}^{\ast
}=X_{\sigma}^{\ast }$ for every stopping time $\sigma$.

For every $\epsilon > 0$ define a stopping time
\begin{equation*}
\mu _{\varepsilon }^{\ast }=\inf \left\{ t\geq 0:V_{t}^{\ast }\leq
a_{t}^{\ast }+\frac{\varepsilon }{35}\right\}.
\end{equation*}
By definition of $\tau $ and $\Gamma ^{\ast }$, one has $\mu _{\varepsilon
}^{\ast }\leq \tau .$

Lepeltier and Mainguenau (1984, Theorem 13) and Proposition \ref{proposition
1} imply that $\mu _{\varepsilon }^{\ast }$ is $\frac{\varepsilon }{35}$%
-optimal for Player 1 in any game $\widetilde{\Gamma }=(\Omega ,\mathcal{A}%
,P,(\mathcal{F}_{t})_{t\geq 0},(a_{t}^{\ast },b_{t}^{\ast },d_{t})_{t\geq
0})$, where $\left( d_{t}\right) _{t\geq 0}$ is any process satisfying $%
a_{t}^{\ast }\leq d_{t}\leq b_{t}^{\ast }$ for any $t\geq 0.$

Let $\delta $ be such that $P(\{|a_{t}-a_{\mu _{\varepsilon }^{\ast }}|<%
\frac{\varepsilon }{35},\ \forall t\in \lbrack \mu _{\varepsilon }^{\ast
},\mu _{\varepsilon }^{\ast }+\delta ]\})>1-\frac{\varepsilon }{35}.$ From
the proof of Proposition \ref{proposition 1} we deduce that $\mu
_{\varepsilon }^{\ast }+r\delta $ is $\frac{\varepsilon }{7}$-optimal for
player 1 in $\Gamma ^{\ast }.$

Now, assume that $\delta $ is sufficiently small so that

\begin{itemize}
\item  $P(\{\mu _{\varepsilon }^{\ast }+\delta <\tau \})\geq P(\{\mu
_{\varepsilon }^{\ast }<\tau \})-\frac{\varepsilon }{7M}$, and

\item  $P(\{|a_{t}-a_{\tau }|<\frac{\varepsilon }{7},|b_{t}-b_{\tau }|<\frac{%
\varepsilon }{7}\ \ \ \forall t\in \lbrack \tau ,\tau +\delta ]\})>1-\frac{%
\varepsilon }{7}$.
\end{itemize}

Define a mixed strategy $\phi _{\varepsilon }$ as follows.
\begin{equation*}
\phi _{\varepsilon }(r,\cdot )=\left\{
\begin{array}{lll}
\mu _{\varepsilon }^{\ast }+r\delta & \ \ \ \ \  & \{\mu _{\varepsilon
}^{\ast }<\tau \}\cup A_{0}, \\
\tau &  & \{\mu _{\varepsilon }^{\ast }=\tau \}\cap \left( A_{1}\cup
A_{2}\right), \\
\mu _{\varepsilon }^{\ast }+r\delta &  & \{\mu _{\varepsilon }^{\ast }=\tau
\}\cap A_{3}.
\end{array}
\right.
\end{equation*}

The proof of Theorem \ref{theorem 1} implies that $\phi _{\varepsilon }$ is $\varepsilon$%
-optimal for player 1 in $\Gamma .$

Assume that $c_\tau \geq b_\tau$ a.s.~(or, equivalently, that $A_1 \cup A_2 =
\Omega$). By the proof of Proposition \ref{proposition 1}, it is $0$-optimal
for Player 1 in $\Gamma $ to stop at time $\tau$, provided the game reaches
time $\tau$. If in addition one has $a_{t}\leq c_{t}\leq b_{t}$ for every $%
t\in \lbrack 0,\tau \lbrack ,$ by the proof of Proposition \ref{proposition
1} and Lepeltier and Mainguenau (1984, Theorem 13) we deduce that the pure
stopping time $\inf \left\{ t \geq 0:V_{t}^{\ast }\leq a_{t}^{\ast
}+\varepsilon \right\} $ is $\varepsilon$-optimal for Player 1 in $\Gamma .$
Hence one obtains the following.

\begin{proposition}
If $c_{t}\in \mathrm{co}\{a_{t},b_{t}\}$ for every $t\in \left[ 0,\tau %
\right] $ then the value exists in pure strategies. An $\varepsilon$-optimal
strategy for Player 1 is $\inf \left\{ t\geq 0:V_{t}^{\ast }\leq a_{t}^{\ast
}+\varepsilon \right\} $, and an $\varepsilon$-optimal strategy for Player 2
is $\inf \left\{ t:V_{t}^{\ast }\geq b_{t}^{\ast }-\varepsilon \right\} .$
\end{proposition}

\begin{corollary}
Every stopping game such that $\left( a_{t},b_{t},c_{t}\right) _{t\geq 0}$
are continuous and satisfies $c_{0}\in \mathrm{co}\{a_{0},b_{0}\}$ admits a
value in pure strategies.
\end{corollary}

\subsection{On final payoff}

Our convention is that the payoff is $0$ if no player ever stops. In fact,
one can add a final payoff as follows. Let $\chi$ be an $\mathcal{A}$%
-measurable and integrable function. The expected payoff that corresponds to
a pair of pure strategies $\left( \mu ,\nu \right) $ is:
\begin{equation*}
{\mathbf{E}}_{P}[a_{\mu }\mathbf{1}_{\{\mu <\nu \}}+b_{\nu }\mathbf{1}%
_{\{\mu >\nu \}}+c_{\mu }\mathbf{1}_{\{\mu =\nu <+\infty \}}+\chi \mathbf{1}%
_{\{\mu =\nu =+\infty \}}].
\end{equation*}

The expected payoff can be written as:
\begin{eqnarray*}
{\mathbf{E}}_{P}\left[ \chi \right] +{\mathbf{E}}_{P}\left[ \left( a_{\mu }-{%
\mathbf{E}}_{P}^{\mathcal{F}_{\mu}} \bigl[ \chi \right] \right) \mathbf{1}%
_{\{\mu <\nu \}}&+&\left( b_{\nu }-{\mathbf{E}}_{P}^{\mathcal{F}_{\nu }}%
\left[ \chi \right] \right) \mathbf{1}_{\{\mu >\nu \}} \\
&+&\left( c_{\mu }-{\mathbf{E}}_{P}^{\mathcal{F}_{\mu }}\left[ \chi \right]
\right) \mathbf{1}_{\{\mu =\nu <+\infty \}}\bigr] ,
\end{eqnarray*}
where ${\mathbf{E}}_{P}^{\mathcal{F}_{\mu }}[\chi ]$ is the conditional
expectation of $\chi $ given the $\sigma $-algebra $\mathcal{F}_{\mu }$.

Define a process $d_{t}:={\mathbf{E}}_{P}^{\mathcal{F}_{t}}\left[ \chi %
\right] .$ Since the filtration satisfies the ``usual conditions'', $\left(
d_{t}\right) _{t\geq 0}$ is a right-continuous martingale (see, e.g.,
Dellacherie and Meyer, 1980, \S VI-4). Hence we are reduced to the study of
the standard stopping game $\Gamma ^{\ast }=(\Omega ,\mathcal{A},P,(\mathcal{%
F}_{t})_{t\geq 0},(a_{t}^{\ast },b_{t}^{\ast },c_{t}^{\ast })_{t\geq 0})$
with $a_{t}^{\ast }=b_{t}-d_{t},$ $b_{t}^{\ast }=b_{t}-d_{t}$ and $%
c_{t}^{\ast }=c_{t}-d_{t}$.

\subsection{On cumulative payoff}

In our definition, players receive no payoff before the game stops. One can
add a cumulative payoff as follows. Let $\left( x_{t}\right) _{t\geq 0}$ be
a progressively measurable process satisfying ${\mathbf{E}}_{P}\left[
\int_{0}^{+\infty }\left| x_{t}\right| dt\right] <+\infty $, and suppose
that the expected payoff that corresponds to a pair of pure strategies $(\mu
,\nu )$ is given by
\begin{equation*}
{\mathbf{E}}_{P}\left[a_{\mu }\mathbf{1}_{\{\mu <\nu \}}+b_{\nu }\mathbf{1}%
_{\{\mu >\nu \}}+c_{\mu }\mathbf{1}_{\{\mu =\nu <+\infty \}}+
\int_{0}^{\min\{\mu, \nu\} }x_{t}dt\right].
\end{equation*}

The expected payoff can be written as
\begin{eqnarray*}
&&{\mathbf{E}}_{P}\left[\left( a_{\mu }+\int_{0}^{\mu }x_{t}dt\right) \mathbf{1%
}_{\{\mu <\nu \}}
+\left( b_{\nu }+\int_{0}^{\nu }x_{t}dt\right) \mathbf{1}%
_{\{\mu >\nu \}}\right.\\
&&\ \ \ \ \ \ +\left.\left( c_{\mu }+\int_{0}^{\mu }x_{t}dt\right) \mathbf{1}%
_{\{\mu =\nu <+\infty \}}\right] +
{\mathbf{E}}_{P}\left[\mathbf{1}_{\{\mu=\nu=+\infty\}} \times \int_0^\infty x_t dt \right].
\end{eqnarray*}
Thus, the game is equivalent to the stopping game $\Gamma ^{\ast }=(\Omega ,%
\mathcal{A},P,(\mathcal{F}_{t})_{t\geq 0},(a_{t}^{\ast },b_{t}^{\ast
},c_{t}^{\ast })_{t\geq 0})$ with terminal payoff $\chi = \int_0^\infty x_t dt $,
where $a_{t}^{\ast }=a_{t}+\int_{0}^t x_s ds$,
$b_{t}^{\ast }=b_{t}+\int_{0}^t x_s ds$, and $c_{t}^{\ast }=c_{t}+\int_{0}^t x_s ds$.

\end{document}